\pdfoutput=1

\documentclass{article}
\usepackage{amssymb, amsmath, mathrsfs, bbm, hyperref}
\usepackage{enumerate}
\usepackage{enumitem}
\usepackage[all, color, pdf]{xy}
\usepackage[all, table]{xcolor}
\usepackage{hyperref}
\usepackage{color}
\usepackage[utf8]{inputenc}
\usepackage[T1]{fontenc} 
\usepackage{mathtools}
\usepackage{listings}
\DeclarePairedDelimiter\abs{\lvert}{\rvert}

\usepackage{ntheorem}

\usepackage[style=alphabetic, backend=bibtex]{biblatex}
\bibliography{MyLibrary}

\theoremstyle{plain}

\newtheorem{axiom}{Axiom}

\newtheorem{example}{Example}

\newtheorem{theorem}{Theorem}
\newtheorem{definition}[theorem]{Definition}
\newtheorem{proposition}[theorem]{Proposition}

\newtheorem{lemma}[theorem]{Lemma}
\newtheorem{observation}[theorem]{Observation}

\newtheorem{remark}{Remark}[section]
\newtheorem{implementation}[remark]{Implementation}

\def\node{\bullet} 
\newcommand{\braket}[1]{\langle #1 \rangle}
\def\glob{{\mathcal G}}
\def\cube{{\mathcal Q}}
\def\O{{\mathcal O}}

\newcommand{\cd}[2][]{\vcenter{\hbox{\xymatrix#1{#2}}}}
\newcommand{\xyg}[2][]{\vcenter{\hbox{\xygraph#1{#2}}}}

\newcommand{\ltwocell  }[3][0.5]{\ar@{}[#2] \ar@{=>}?(#1)+/r 0.2cm/;?(#1)+/l 0.2cm/^{#3}}
\newcommand{\rtwocell  }[3][0.5]{\ar@{}[#2] \ar@{=>}?(#1)+/l 0.2cm/;?(#1)+/r 0.2cm/^{#3}}
\newcommand{\dtwocell  }[3][0.5]{\ar@{}[#2] \ar@{=>}?(#1)+/u 0.2cm/;?(#1)+/d 0.2cm/^{#3}}
\newcommand{\utwocell  }[3][0.5]{\ar@{}[#2] \ar@{=>}?(#1)+/d 0.2cm/;?(#1)+/u 0.2cm/_{#3}}
\newcommand{\congcell  }[2][0.5]{\ar@{}[#2]|(#1){\cong}}
\newcommand{\drtwocell  }[3][0.5]{\ar@{}[#2] \ar@{=>}?(#1)+/u 0.14cm/+/l 0.14cm/;?(#1)+/d 0.14cm/+/r 0.14cm/^{#3}}

\newcommand\overlap[2]{\mathrel{\ooalign{\hss$#1$\hss\cr$#2$}}}

\lstdefinestyle{customc}{
backgroundcolor=\color{white},
breakatwhitespace=true,
breaklines=true,
keepspaces=true,
belowcaptionskip=1\baselineskip,
breaklines=true,
xleftmargin=\parindent,
showstringspaces=false,
basicstyle=\footnotesize\ttfamily,
keywordstyle=\bfseries,
}
\lstnewenvironment{code}[1][]%
  {\\\minipage{0.9\linewidth}}
  {\endminipage\\}

\title{Formalizing parity complexes}

\author{Mitchell Buckley}

\xyoption{rotate}

\begin{document}

\maketitle

\begin{abstract}
We formalise, in Coq, the opening sections of \emph{Parity Complexes}~\cite{Street1991} up to and including the all important excision of extremals algorithm.
Parity complexes describe the essential combinatorial structure exhibited by simplexes, cubes and globes, that enable the construction of free $\omega$-categories on such objects.
The excision of extremals is a recursive algorithm that presents every cell in such a category as a (unique) composite of atomic cells. 
This is the sense in which the $\omega$-category is (freely) generated from its atoms.
Due to the complicated multi-dimensional nature of this work, the detail of definitions and proofs can be hard to follow and verify.
Indeed, some corrections were required some years following the original publication~\cite{Street1994}.
Our formalisation verifies that all cases of each result operate as stated.
In particular, we indicate which portions of the theory can be proved directly from definitions, and which require more subtle and complex arguments.
By identifying results that require the most complicated proofs, we are able to investigate where this theory might benefit from further study and which results need to be considered most carefully in future work.
\end{abstract}


\section{Introduction} \label{sec:intro}

An $n$-simplex $\Delta_n$ is a geometric figure that generalises the notion of triangle or tetrahedron to $n$-dimensional space.
Simplexes have a number of properties that make them useful in algebraic topology, algebraic geometry and homotopy theory where they often play a foundational role.
Each $n$-simplex can be oriented in such a way that it forms an $n$-category.
We include below the cases for $n=1$, $2$ and $3$.
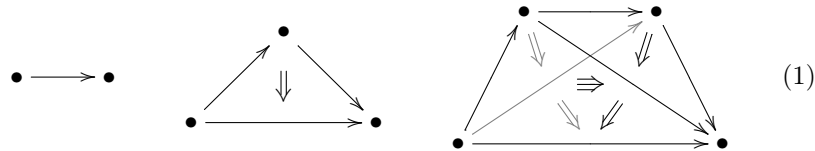
\begin{equation}
\qquad
\cd[]{
\node \ar[r] & \node
}
\qquad 
\cd[]{
& \node \ar[dr] \dtwocell[0.6]{d}{}  &  \\
\node \ar[ur] \ar[rr] & & \node
}
\qquad 
\xyg{
!{<0em,0em>;<2.5em,0em>:<0em,2.5em>::}
{\node}="p0" [uur]
{\node}="p1" [rr]
{\node}="p2" [ddr]
{\node}="p3"
"p0":"p1"|-{}="p01"
"p0":@[gray]"p2"|(0.45){}="p02"
"p0":"p3"|-{}="p03"
"p1":"p2"|-{}="p12"
"p1":"p3"|(0.55){}="p13"
"p2":"p3"|-{}="p23"
"p1":@{}"p02"|(0.3){}="p012a"|(0.75){}="p012b"
"p2":@{}"p13"|(0.3){}="p123a"|(0.75){}="p123b"
"p13":@{}"p03"|(0.3){}="p013a"|(0.75){}="p013b"
"p02":@{}"p03"|(0.3){}="p023a"|(0.75){}="p023b"
"p012a":@{=>}@[gray]"p012b"
"p023a":@{=>}@[gray]"p023b"
"p013a":@{=>}"p013b"
"p123a":@{=>}"p123b"
"p13":@{}"p02"|(0.35){}="p0123a"|(0.65){}="p0123b"
"p0123b":@3{->}"p0123a"
}
\end{equation}
At low dimensions, it is not hard make each of these into an $n$-category. 
At higher dimensions, say $n > 3$, it is quite hard to describe the $n$-category structure because the source and target of each cell are large pasting diagrams in high dimensions.

Beginning in the late 1970's Ross Street, together with John Roberts and Jack Duskin, began investigating how this process could be rigorously extended to any $n$.
This was achieved in~\cite{Street1987} where the process was described for the simplexes and the corresponding categories were dubbed the \emph{orientals} (referring to the fact that they are oriented).
The main motivation at this time stemmed from non-abelian cohomology where various constructions rely on the orientals.

At the same time, Iain Aitcheson was developing a similar series of results for $n$-cubes: that each cube could be given an orientation in such a way that it forms an $n$-category or even an $\omega$-category~\cite{Aitchison1986}.
A third example of this phenomenon is found in $n$-globes where the corresponding $n$-categories have a very simple description.
For more on the usefulness of simplexes and cubes, see Street's survey~\cite{Street1995}.

Following these successes, the goal was then to describe the general structure of all oriented multi-dimensional structures for which it is possible to extract free $\omega$-categories in the style of these three examples.
The early 1990's yielded a number of related solutions.
Ross Street defined a structure called a \emph{parity complex} and gave an explicit description of the $\omega$-category associated to each parity complex~\cite{Street1991}. 
Some minor corrections were added in~\cite{Street1994}.
Richard Steiner contributed \emph{directed complexes} as a generalisation of directed graph.
He showed that loop-free directed complexes generated free $\omega$-categories in the appropriate way~\cite{Steiner1993}. 
Both of these authors also showed that their respective structures were closed under product and join and covered the three main examples of simplexes, cubes, and globes.
Around the same period Mike Johnson was working on a formal description of \emph{pasting scheme} for  $\omega$-categories~\cite{Johnson1988,Johnson1989}, and was able to describe the free $\omega$-category on such structures. 
He included the simplexes as his primary example, and there is a strong sense in which this addressed the same problem.
Further related work can be found in~\cite{Al-Agl1993,Steiner2004} and also in~\cite{Verity2008} where a conjecture of Street--Roberts is proved in the closing chapter.

Our interest centres on \emph{Parity Complexes} which takes a particularly `hands-on' approach and describes the combinatorics of this construction in full detail.
Our goal is to encode and verify the opening sections of the article up to the excision of extremals algorithm.
The theory shows how to build, for any parity complex $C$ an $\omega$-category $\O(C)$.
The excision of extremals algorithm shows that each cell can be presented as a (unique) composite of atomic cells; this is the sense in which $\O(C)$ is (freely) generated from the atoms. 
The algorithm can also be used to generate explicit algebraic descriptions of the cells in $\O(C)$.

Our motivation is two-fold.
First, some of the combinatorial arguments in Street's text can be difficult to follow and can easily conceal errors; this is illustrated by the fact that corrections were later required. 
We will provide some confirmation that the corrections have addressed all issues.
Second, a computer-verified encoding provides a good resource for understanding the intricacies of these complicated structures and opens a path to further refinements of the material.
We have not attempted to formalized the entirety of the theory.
The essential combinatorics are contained in sections 1 to 4 and culminate in the excision of extremals algorithm which is the final result that we encode.

From this point on we often refer to~\cite{Street1991} as `the original text', and to~\cite{Street1994} as `the corrigenda'.
The content extracted directly from the original article is numbered consecutively with the original name included parentheses.
Sections within this article are referred to as `Section $n$' and those of the original text are referred to as `\S$n$' or `Section $n$ of~\cite{Street1991}'.

We programmed everything in Coq~\cite{coq} and the code is freely available for inspection at the following location.
\begin{quote}
\url{https://github.com/MitchellBuckley/Parity-Complexes}
\end{quote}

In Section~\ref{sec:foundations} we outline the foundational mathematics that needs to be introduced for an encoding of parity complexes. 
We also outline how we chose to implement this foundation.
In Sections~\ref{sec:definitions} to~\ref{sec:product} we outline the main combinatorial content of~\cite{Street1991} section by section. 
At each stage we comment on the intuition underlying each result and discuss our implementation of the definitions and results. 
We pay particular attention to those parts of the material that were difficult to translate into Coq.
Though our encoding focuses on \S1 to \S4, we do comment on some material in \S5.
In Section~\ref{sec:lessons} we outline the few lessons we have learned in computer-verified encoding of mathematics.
Section~\ref{sec:conclusion} contains concluding remarks.

\section{Required Foundations} \label{sec:foundations}

Parity complexes are described using basic set theory and partially ordered sets.
In particular, we must implement: 
\begin{itemize}
\item sets;
\item set union, set intersection, set difference, etc.;
\item finite sets;
\item cardinality of finite sets;
\item partial orders; and
\item segments of partial orders.
\end{itemize}
Many of these structures are already encoded in different parts of the Coq standard library but for various reasons we have chosen to reimplement much of that material.

\subsection{Sets}
We implemented sets using the \lstinline{Ensembles} standard library.
This involves a universe type \lstinline{U : Type} on which all our sets will be based.
Then a set is an \emph{Ensemble}: an indexed proposition \lstinline{U $\rightarrow$ Prop}.
An element of the universe \lstinline{x : U} is a member of a set \lstinline{A : U $\rightarrow$ Prop} when the corresponding proposition \lstinline{A x} is true.
Inclusion of sets relies on logical implication.
\begin{code}
Definition Ensemble := U $\rightarrow$ Prop.
Definition In (A:Ensemble) (x:U) : Prop := A x.
Definition Included (B C:Ensemble) : Prop := 
  forall x:U, In B x $\rightarrow$ In C x.
\end{code}
Set operations union, intersection, and set difference are all implemented using point-wise logical operations: 
\begin{code}
Union A B        := fun x => (A x $\vee$ B x)
Intersection A B := fun x => (A x $\wedge$ B x)
Setminus A B     := fun x => (A x $\wedge$ $\neg$(B x))
\end{code}
For the purposes of this section we suppose that we always work with a fixed universe \lstinline{U}.

The Coq language has a convenient feature that allows us to introduce notation for these operations.
\begin{code}
Notation "x $\in$ B"       := (In A x) (at level 71).
Notation "A $\subseteq$ B"       := (Included A B) (at level 71).
Notation "A $\cup$ B"       := (Union A B) (at level 61).
Notation "A $\cap$ B"       := (Intersection A B) (at level 61).
Notation "A '\' B"     := (Setminus A B) (at level 61).
\end{code}
Each special symbol is introduced as a utf-8 character which Coq has no problem recognising.
This feature makes the code much more readable.

\subsection{Finiteness and cardinality}
Finiteness is implemented using the same basic idea as the \lstinline{Finite_sets} standard library.
This library contains an inductively defined proposition \lstinline{Finite} stating that a set $S$ is finite when $S = \emptyset$, or $S = \{x\} \cup S'$ where $S'$ is finite and $x \not\in S'$.
Cardinality is implemented in a similar way: there is an inductively defined proposition \lstinline{cardinal} stating that a set $S$ has cardinality $0$ when it is empty and has cardinality $n+1$ when $S = \{ x \} \cup S'$, $x \not\in S'$, and $S'$ has cardinality $n$.
\label{code:finite}
\begin{code}
Inductive Finite : Ensemble U $\rightarrow$ Prop :=
  | Empty_is_finite : Finite (Empty_set U)
  | Union_is_finite :
      forall A:Ensemble U,
        Finite A $\rightarrow$ forall x:U, $\neg$ In U A x $\rightarrow$ 
          Finite (Add U A x).
          
Inductive cardinal : Ensemble U $\rightarrow$ nat $\rightarrow$ Prop :=
  | card_empty : cardinal (Empty_set U) 0
  | card_add :
      forall (A:Ensemble U) (n:nat),
        cardinal A n $\rightarrow$ forall x:U, $\neg$ In U A x $\rightarrow$ 
          cardinal (Add U A x) (S n).
\end{code}

Cardinality and finiteness are related by the result \lstinline{forall S, (Finite S <-> exists n, cardinal S n)}.
When our universe has decidable equality we can show that finiteness interacts well with set operations, for example \lstinline{forall A B, Finite A $\wedge$ Finite B $\rightarrow$ Finite (A $\cup$ B)}.
When our universe does not have decidable equality, there are some classical results that may not hold for all finite sets. 
For example $\neg\neg S = S$ and $S \subseteq T \to T = (T \setminus S) \cup S$. 
We don't know exactly how the theory will change if this avenue is explored. 
Since most types of interest have decidable equality, we are not worried about including this as an assumption in our formalisation.

Note that the definitions we have given here are those contained in the standard library while our implementation is slightly different; the differences are explained in Section \ref{equality_of_sets} below.

\subsection{Partial orders}
Some material on partial orders is available in the \lstinline{Relations} standard library.
Our particular requirements for orders were slightly more complicated than that library allowed for and we found it simpler to explicitly prove basic results as they were needed.

\subsection{Equality of sets}\label{equality_of_sets}
We say that two sets $S$ and $T$ are \emph{equal} when they are equal as terms of the type \lstinline{Ensemble U}; that is, they are equal as indexed propositions and \lstinline{forall x, S x = T x}.
We write $S = T$ to indicate that $S$ and $T$ are equal.
This is the standard notion that is built into Coq and allows us to replace $S$ with $T$ in any expression.

There is another notion of equality: we say that $S$ and $T$ are \emph{the same} when they contain the same elements.
This is the usual notion of set equality used in mathematics.
Equivalently, two sets are the same when they are equivalent as indexed propositions (\lstinline{forall x, S x $\leftrightarrow$ T x}), or when \lstinline{S $\subseteq$ T} and \lstinline{T $\subseteq$ S}.
We write \lstinline{Same_set S T} or \lstinline{S == T} to indicate that $S$ and $T$ are the same.

If two sets are equal then they are certainly the same but two sets can be the same without being equal.
For example, the sets \lstinline{fun x => x = 0} and \lstinline{fun x => 1 + x = x} in \lstinline{Ensemble nat} are the same but not equal.
The standard library \lstinline{Ensembles} contains an extensionality axiom stating that 
\lstinline{forall A B, A == B $\rightarrow$ A = B}.
In order to keep our formalisation as constructive as possible we are careful never to use the axiom in our formalisation.

The standard facilities of Coq will allow for rewriting $S$ with $T$ whenever $S$ is equal to $T$.
However, when $S$ is the same as $T$ we do not have any guarantee that such rewrites are legitimate.
In this situation, our type of sets \lstinline{Ensemble U} becomes a \emph{setoid}: a set equipped with an equivalence relation.
Then \lstinline{S == T} implies that we may rewrite \lstinline{S} for \lstinline{T} in any expression which is built up of operations that preserve the equivalence relation.
This rewrite facility is provided by the standard library \lstinline{Setoid} and requires us to prove that \lstinline{Same_set} is an equivalence relation and that the appropriate set operations preserve the equivalence.

Without the extensionality axiom it is not possible to prove that \lstinline{Finite S} and \lstinline{S == T} implies \lstinline{Finite T}.
Consider the usual inductive definition of \lstinline{Finite} given on page~\ref{code:finite}. 
By examining the constructors we can see that the proposition holds for sets of the form \lstinline{Emptyset}, \lstinline{Add Emptyset x}, and \lstinline{Add (Add Emptyset x) y} etc.
It is easy to describe sets that are finite, but not equal to sets of this kind.
Take \lstinline{fun n => (n < 2)}, or \lstinline{fun n => (S n = 0)}; while these sets are \emph{the same} as  \lstinline{Add (Add Emptyset 0) 1} and \lstinline{Emptyset} they are not equal to those sets, and thus not finite under the usual definition.
Something similar happens with the definition of finite cardinality.
This problem can be solved in more than one way.
We chose to solve this by adding a third constructor for \lstinline{Finite} that explicitly introduces the property that \lstinline{Finite S $\wedge$ S == T $\rightarrow$ Finite T}.
This modification allows us to recover this basic property of finite sets without the extensionality axiom.

\subsection{More on finiteness}
In many cases we augmented the standard library with extra results about finite sets that were not already present.
We found that setting up this basic theory was often tedious, but occasionally an enjoyable exercise in constructive mathematics.
For instance, it became clear at some point that certain basic results about sets could not be proved without supposing that equality in \lstinline{U} is decidable, i.e.\ \lstinline{forall (a b : U), (a=b) $\vee$ $\neg$(a=b)}.
Since none of the examples used here or in the literature need a universe \lstinline{U} without decidable equality, we have made this a further assumption in our implementation.

If one wanted to reason about, say sets of integer sequences, then the obvious universes to use would be $nat -> Int$ or $Stream Int$ each of which lacks decidable equality. 
In that case one would find that various simple results concerning finite sets would not hold.

\subsection*{} 
We have now covered the essential mathematical foundations required for a formalisation of parity complexes.
More details can be found by examining the code itself.
In the following three sections we summarise \S1 to \S4 of \emph{Parity Complexes} together with modifications given in the corrigenda.
This content is sufficient to express the excision of extremals algorithm (Theorem~\ref{theorem_4_1}).
As we progress through the material we will reproduce definitions and terminology almost verbatim from~\cite{Street1991,Street1994}.
In each case we will explain the underlying intuition of the material, comment on our implementation, and indicate where our formalisation shed light on the underlying arguments.
We also discuss some aspects of the final sections (5 and 6) of~\cite{Street1991} though none of that material was been formalised.


\section{Definitions and the simplex example}  \label{sec:definitions}

We begin by summarising the content of Section 1 of~\cite{Street1991}.

\begin{definition}
A \emph{parity complex} is a graded set
\begin{equation}
C = \sum_{n=0}^{\infty} C_n
\end{equation}
together with, for each $x \in C_{n+1}$ two disjoint, non-empty, finite sets ${x}^{+}, {x}^{-} \subseteq C_n$ subject to Axioms 1, 2, 3A and 3B which appear below.
\end{definition}

From this point onward we will work exclusively within a single parity complex $C$ as described above. 
When we say $S \subseteq C$ we mean that $S$ is a subset of the underlying graded set of the parity complex.
When we say $x \in C$ we mean that $x$ is an element of the underlying graded set of the parity complex.

Before we list the axioms we will introduce some terminology.
If $x \in C$ then elements of ${x}^{-}$ are called \emph{negative faces} of $x$, and those of ${x}^{+}$ are called \emph{positive faces} of $x$.
We will sometimes refer to $x^-$ and $x^+$ as \emph{face-sets} of $x$.
Given $S \subseteq C$, let $S^-$ denote the set of elements of $C$ which occur as negative faces of some $x\in S$, and similarly for $S^+$; thus
\begin{equation}
S^- = \bigcup\limits_{w \in S}^{} w^- \quad \text{and} \quad S^+ = \bigcup\limits_{w \in S}^{} w^+\rlap{\ .}
\end{equation}
Each subset $S \subseteq C$ is graded via $S_n = S \cap C_n$. 
The \emph{$n$-skeleton} of $S \subseteq C$ is defined by
\begin{equation}
S^n := \sum^{n}_{k=0} S_k \rlap{\rlap{~.}}
\end{equation}
Call $S$ \emph{$n$-dimensional} when it is equal to its $n$-skeleton.

The broad intuition is to see this structure as a generalisation of directed graph.
Elements of $C_0$ are vertices, elements of $C_1$ are directed edges, elements of $C_2$ are directed `faces', elements of $C_3$ are directed `volumes', and so on.
The usual notion of source and target are replaced by face-sets $x^-$ and $x^+$.
The following is a basic example of this structure of dimension two.
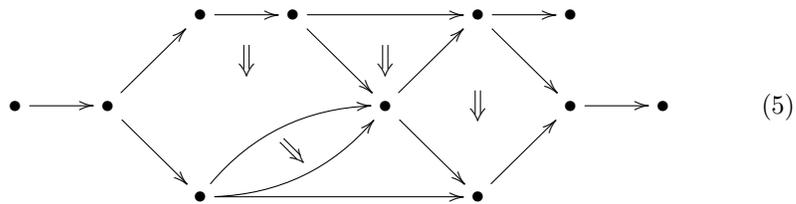
\begin{equation}\label{general_example}
\cd[]{
 & & 
 \node \ar[r] & 
 \node \ar[rr] \ar[dr] & & 
 \node \ar[r] \ar[dr] & 
 \node \\
 \node \ar[r] & 
 \node \ar[ur] \ar[dr] \dtwocell{urrr}{} & & \dtwocell{urr}{} & 
 \node \ar[ur] \ar[dr] \dtwocell{rr}{} & & 
 \node \ar[r] & 
 \node \\   
 & & 
 \node \ar@/^1em/[urr] \ar@/_1em/[urr] \ar[rrr] \drtwocell{urr}{} & & & 
 \node \ar[ur] &  
}
\end{equation}
Notice that elements above dimension 1 can have more than one source-face or target-face.

Without the axioms below, this structure is very general indeed and many unusual examples can be provided. 
When the axioms are applied, possible examples become much better behaved.
Examples of arbitrary dimension can be constructed from simplexes, cubes, and other kinds of polytopes as seen below.
Of course, the simplexes provide the main motivation for understanding these kinds of structures.

So far we have described the data of a parity complex: a graded set with a pair of face-set maps $(-)^-,(-)^+ \colon C_{n+1} \to \mathcal{P}(C_n)$.
We now describe the required axioms.

\begin{axiom}\label{axiom1}
For all $x\in C$,
$${x}^{++} \cup {x}^{--} = {x}^{-+} \cup {x}^{+-}$$
where $x^{++} = (x^+)^+$ etc.
\end{axiom}
This is a kind of globularity condition that ensures various face-sets are appropriately related. 
The following diagram is an example where $x \in C_2$ and both $x^-$ and $x^+$ have four elements.
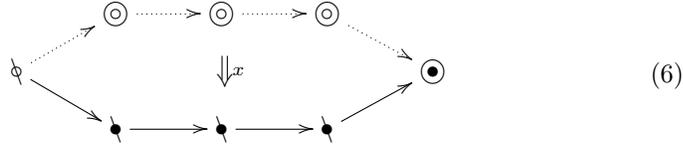
\begin{equation}\label{axiom1_diagram}
\cd[@R-1.5em]{ & 
	\overlap{\circ}{\bigcirc} \ar@{.>}[r] & 
	\overlap{\circ}{\bigcirc} \ar@{.>}[r] & 
	\overlap{\circ}{\bigcirc} \ar@{.>}[dr] &  \\
	\overlap{\circ}{\setminus} \ar@{.>}[ur] \ar[dr] & & & & 
	\dtwocell{llll}{x} \overlap{\bullet}{\bigcirc} \\ & 
	\overlap{\bullet}{\setminus} \ar[r] & 
	\overlap{\bullet}{\setminus} \ar[r] & 
	\overlap{\bullet}{\setminus} \ar[ur] &   	 
}
\end{equation}
Edges marked with a dotted line belong to $x^-$, the other edges belong to $x^+$.
Vertices marked with a $\bullet$ belong to $x^{++}$, those marked with a $\circ$ belong to $x^{--}$, those marked with a $\bigcirc$ belong to $x^{-+}$, and those marked with a $\setminus$ belong to $x^{+-}$.
In particular, this axiom implies that ${x}^{++} \subseteq {x}^{-+} \cup {x}^{+-}$, that is, positive faces of positive faces must be the negative face of a positive face, or the positive face of a negative face.

Notice that in both~\eqref{general_example} and~\eqref{axiom1_diagram} the set of source (target) faces have all elements aligned in a common direction and they do not branch apart.
This behaviour is guaranteed by introducing Axiom~\ref{axiom2} below.

Suppose that $S$ and $T$ are subsets of $C$. 
We write $S \perp T$ when $S^-\cap T^- = S^+\cap T^+ = \emptyset$. 
This extends to elements by $x \perp y$ when $x^-\cap y^- = x^+\cap y^+ = \emptyset$\ \footnote{This could equivalently be defined for elements first and then extended to sets afterwards.}. 
A subset $S \subseteq C$ is called \emph{well-formed} when $S_0$ has at most one element, and, for all $x,y\in S_n$ $(n > 0)$, if $x \not= y$ then $x\perp y$.
Broadly speaking, a set is well-formed when its elements do not form any branchings like
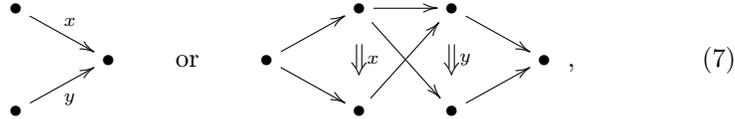
\begin{equation}\label{diag:branchings}
\cd[@R-1.5em]{
	\node \ar[dr]^{x} &    \\
	 & \node \\
	\node \ar[ur]_{y} & 
}
\qquad \text{or} \qquad
\cd[@R-1.5em]{
	& \node \ar[r] \dtwocell{dd}{x} \ar[ddr] & \node \ar[dr] \dtwocell{dd}{y} &   \\
\node \ar[ur] \ar[dr] & & & \node \\
	& \node \ar[uur] & \node \ar[ur] &   
}
\rlap{\ ,}
\end{equation}
and it contains at most one element of dimension 0.
In each of the diagrams above we can observe that $\{x,y\}$ is not well-formed, while $\{x\}$, $\{y\}$, $x^+$, $x^-$, $y^+$, and $y^-$ are all well-formed.
This diagram depicts branchings in dimensions 1 and 2, but well-formedness prevents branching in all dimensions.
The condition on dimension zero does not force parity complexes to have a single element of dimension zero, but that (using the axiom below) elements of dimension $1$ have a single source vertex and a single target vertex.

\begin{axiom}\label{axiom2}
For all $x \in C$, ${x}^{-}$ and ${x}^{+}$ are well-formed.
\end{axiom}
If we think of the union ${x}^{-} \cup {x}^{+}$ as forming a boundary of $x$ as in~\eqref{axiom1_diagram} above then this axiom ensures that the boundary looks something like the boundary of a polytope.
For those familiar with higher categories, this condition ensures that the face-sets look like valid pasting diagrams.

Suppose that $x,y\in C$.
We write $x < y$ whenever ${x}^{+} \cap {y}^{-}$ is non-empty.
That is, when $x$ and $y$ abut by having a common element in their respective sets of positive and negative faces.
This implies $x \not= y$ since $x^-$ and $x^+$ are always disjoint.
We then let $\lhd$ be the reflexive transitive closure of $<$.
An example is 
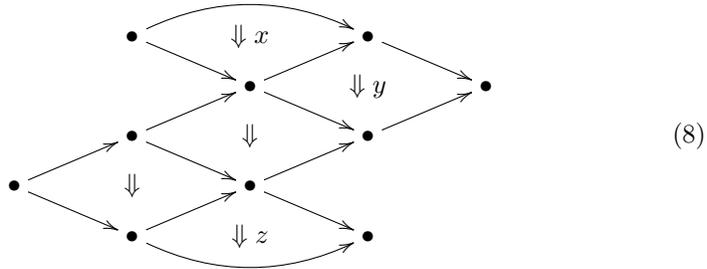
\begin{equation}
\cd[@!C@R-2em]{ & 
\node \ar@/^{1.2em}/[rr] \ar[dr] & 
\Downarrow x & 
\node \ar[dr] &  \\ & & 
\node \ar[ur] \ar[dr] & 
\Downarrow y & 
\node  \\ & 
\node \ar[dr] \ar[ur] & 
\Downarrow & 
\node \ar[ur] & \\
\node \ar[ur] \ar[dr] & 
\Downarrow & 
\node \ar[ur] \ar[dr] & &  \\ & 
\node \ar@/_{1.2em}/[rr] \ar[ur] & 
\Downarrow z & 
\node &  	 
}
\end{equation}
where $x < y$ and $y \lhd z$.
In this case we often say that there is a path from $x$ to $z$.
For all $S \subseteq C$ we let $\lhd_S$ denote the reflexive transitive closure of $<$ restricted to $S$.
When $x \lhd_S z$ we often say there is a path from $x$ to $z$ in $S$.

While Axioms 1 and 2 can be seen as imposing some of the basic structural behaviour of graphs, the following axiom restricts us to certain `loop-free' graphs.
\begin{axiom}\label{axiom3}
For all $x,y\in C$,
\begin{enumerate}[label=\Alph*., itemindent=3.5em, labelsep=1.5em]
\item $x \lhd y \lhd x$ implies $x = y$.
\item if $x\lhd y$ then
$\forall z\in C$, 
$\neg (x \in {z}^{+} \wedge y \in {z}^{-})$ and
$\neg (y \in {z}^{+} \wedge x \in {z}^{-})$.
\end{enumerate}
\end{axiom}
Axiom~\ref{axiom3}A says that $\lhd$ is anti-symmetric, or, that there are no paths that loop within a fixed dimension.
Axiom~\ref{axiom3}B says that there are no paths that cross between the face-sets of any element $z$.
That is, we avoid circumstances where a path can cross from one face-set to the other face-set of an element $z$ as in the diagram below.
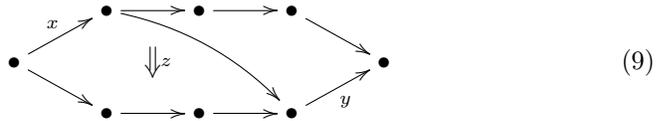
\begin{equation}
\cd[@R-1.5em]{
	& \node \ar[r] \ar@/^.7em/[ddrr] & \node \ar[r] & \node \ar[dr] & \\
	\node \ar[ur]^x \ar[dr] & & &  \dtwocell{lll}{z} & \node \\
	& \node \ar[r] & \node \ar[r] & \node \ar[ur]_y &   	 
}
\end{equation}

These are all the axioms for a parity complex.
The following examples come from p.318--319 of~\cite{Street1991}.

\begin{example}
A 1-dimensional parity complex is precisely a directed graph with no circuits.
\end{example}

\begin{example}
The \emph{$\omega$-glob} is the parity complex $\glob$ defined by 
$\glob_n = \{(\epsilon, n) \colon \epsilon = \ominus \text{ or } \oplus\}$, and 
$ (\epsilon, n+1)^- = \{(\ominus,n)\}$ and $(\epsilon, n+1)^+ = \{(\oplus,n)\}$.
Elements of dimension 0, 1, and 2 are `$n$-discs'. 
There are precisely two elements at each dimension, each of which has exactly one source face and exactly one target face.
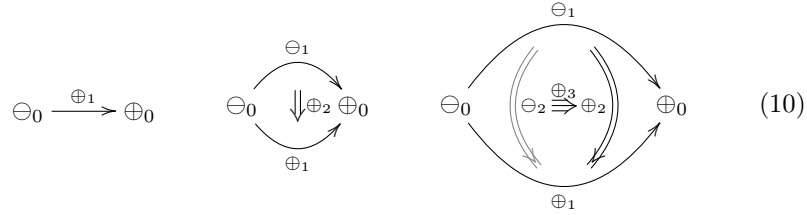
\begin{equation}
\cd[]{
\ominus_0 \ar[r]^{\oplus_1} & \oplus_0
}
\qquad 
\cd[@R-1.5em]{
& \\
\ominus_0 \ar@/^1.6em/[r]^{\ominus_1} \ar@/_1.6em/[r]_{\oplus_1} \dtwocell{r}{\oplus_2} & \oplus_0 \\
&
}
\qquad
\xyg{
!{<0cm,0cm>;<0.8em,0cm>:<0cm,0.8em>::}
!{(0,0)  }*+{\ominus_0}="a0"
!{(10,0) }*+{\oplus_0}="b0"
!{(4, 3) }*+{ }="a2a"
!{(4,-3) }*+{ }="a2b"
!{(6, 3) }*+{ }="b2a"
!{(6,-3) }*+{ }="b2b"
"a0":@/^3em/"b0"^{\ominus_1}
"a0":@/_3em/"b0"_{\oplus_1}
"a2a":@/_1.1em/@{=>}@[gray]"a2b"^{\ominus_2}="a2"
"b2a":@/^1.1em/@{=>}"b2b"_{\oplus_2}="b2"
"a2":@{}"b2"|(0.3){}="a3a"|(0.7){}="a3b"
"a3a":@3{->}"a3b"^{\oplus_3}="a3"
}
\end{equation}
We use $\ominus_n$ and $\oplus_n$ as short-hand for $(\ominus,n)$ and $(\oplus,n)$.
\end{example}

\begin{example}
The \emph{$\omega$-simplex} is the parity complex $\Delta$ described as follows. 
Let $\Delta_n$ denote the set of $(n+1)$-element subsets of the set of natural
numbers $N = \{0,1,2,\dots\}$.
Each $x \in \Delta_n$ is written as $(x_0,x_1,\dots,x_n)$ where $x_0 < x_1 < \dots < x_n$. 
Let $x\delta_i$ denote the set obtained from $x$ by deleting $x_i$.
Take $x^-$ to be $\{x\delta_i : i \text{ odd}\}$ and $x^+$ to be $\{x\delta_i : i \text{ even}\}$.
Elements of dimension 0, 1, and 2 are `$n$-simplexes'.
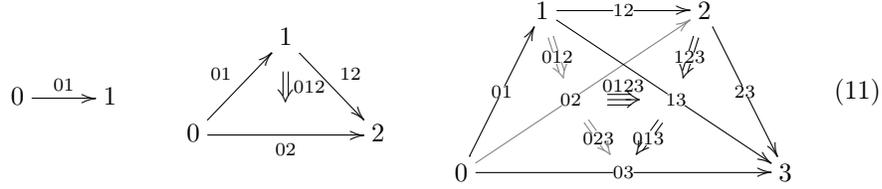
\begin{equation}
\cd[]{
0 \ar[r]^{01} & 1
}
\qquad 
\cd[]{
& 1 \ar[dr]^{12} \dtwocell{d}{012}  &  \\
0 \ar[ur]^{01} \ar[rr]_{02} & & 2 
}
\qquad
\xyg{
!{<0em,0em>;<3em,0em>:<0em,3em>::}
{0}="p0" [uur]
{1}="p1" [rr]
{2}="p2" [ddr]
{3}="p3"
"p0":"p1"|-{01}="p01"
"p0":@[gray]"p2"|(0.45){02}="p02"
"p0":"p3"|-{03}="p03"
"p1":"p2"|-{12}="p12"
"p1":"p3"|(0.55){13}="p13"
"p2":"p3"|-{23}="p23"
"p1":@{}"p02"|(0.3){}="p012a"|(0.75){}="p012b"
"p2":@{}"p13"|(0.3){}="p123a"|(0.75){}="p123b"
"p13":@{}"p03"|(0.3){}="p013a"|(0.75){}="p013b"
"p02":@{}"p03"|(0.3){}="p023a"|(0.75){}="p023b"
"p012a":@{=>}@[gray]"p012b"|-{012}
"p023a":@{=>}@[gray]"p023b"|-{023}
"p013a":@{=>}"p013b"|-{013}
"p123a":@{=>}"p123b"|-{123}
"p13":@{}"p02"|(0.35){}="p0123a"|(0.65){}="p0123b"
"p0123b":@3{->}"p0123a"^{0123}
}
\end{equation}
We use $abcd$ as short-hand for $(a,b,c,d)$ and similarly at other dimensions.
\end{example}

\begin{example}
The \emph{$\omega$-cube} is the parity complex $\cube$ described as follows.
The elements are infinite sequences of the three symbols $\ominus, \odot, \oplus$ containing a finite number of $\odot$'s and ending with an infinite string of $\ominus$'s. 
The dimension of an element is the number of $\odot$'s appearing in it. 
Let $x\delta^-_i$ denote the sequence obtained from $x$ by replacing the $i$-th $\odot$ by $\ominus$ when $i$ is odd and by $\oplus$ when $i$ is even. 
Similarly, $x\delta^-_i$ is defined by interchanging $\ominus$ and $\oplus$ in the previous sentence. 
For $x \in \cube_n$, define $x^\epsilon = \{x\delta_i^\epsilon : 1 < i < n\}$.
The \emph{$n$-cube} is the parity complex built the same way but using only lists of length $n$ (and that are not required to end with an infinite string of $\ominus$'s).
The $n$-cubes of dimension 1, 2, and 3 are displayed below.
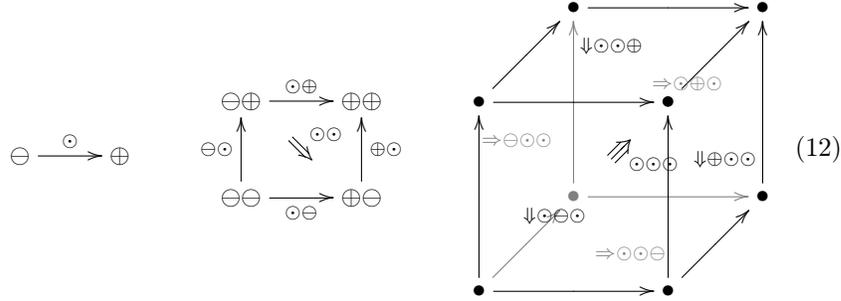
\begin{equation}
\qquad
\cd[]{
\ominus \ar[r]^{\odot} & \oplus
}
\qquad 
\cd[]{
\drtwocell{dr}{\odot\odot} \ominus\oplus \ar[r]^{\odot\oplus} & \oplus\oplus   \\
\ominus\ominus \ar[u]^{\ominus\odot} \ar[r]_{\odot\ominus} & \oplus\ominus \ar[u]_{\oplus\odot}
}
\qquad
\xyg{
!{<0em,0em>;<3.5em,0em>:<0em,3.5em>::}
{\node}="100" [u(1)r(1)]
{\node}="101" [r(2)]
{\node}="111" [d(1)l(1)]
{\node}="110" [d(2)]
{\node}="010" [u(1)r(1)]
{\node}="011" [l(2)]
{\color{gray}\node}="001" [d(1)l(1)]
{\node}="000" [ur]
"001":@[gray]"101"|-{}="a01"
"001":@[gray]"011"|-{}="0a1"
"000":@[gray]"001"|-{}="00a"
"000":"100"|-{}="a00"
"010":"110"|-{}="a10"
"011":"111"|-{}="a11"
"000":"010"|-{}="0a0"
"100":"110"|-{}="1a0"
"101":"111"|-{}="1a1"
"100":"101"|-{}="10a"
"010":"011"|-{}="01a"
"110":"111"|-{}="11a"
%
"001":@{}"110"|(0.4){}="A"|(0.6){}="B"
"100":@{}"001"|(0.4){\color{gray}\Rightarrow \ominus\odot\odot} 
"010":@{}"001"|(0.4){\color{gray}\Rightarrow \odot\odot\ominus} 
"111":@{}"001"|(0.4){\color{gray}\Rightarrow \odot\oplus\odot} 
"000":@{}"110"|(0.4){\Downarrow \odot\ominus\odot} 
"011":@{}"110"|(0.4){\Downarrow \oplus\odot\odot} 
"101":@{}"110"|(0.4){\Downarrow \odot\odot\oplus} 
%
%
"A":@3{->}"B"_{\odot\odot\odot}:@{}"110"\node
}
\end{equation}
Some labels have been omitted from the last diagram in order to keep it readable.
\end{example}

\begin{remark}
It might seem unusual to insist that elements of the $\omega$-cube end with an infinite string of $\ominus$'s.
If we omit this condition then the parity complex would not contain any (finite) paths from $\ominus\ominus\dots$ to $\oplus\oplus\dots$ and so it would become somewhat disconnected.
The underlying type for our implementation would also become undecidable.
We have not investigated the implications of this distinction.
\end{remark}

Before continuing our exposition of \S1 we will comment briefly on our implementation.

\begin{implementation}
The basic data for a parity complex without the axioms is sometimes called a \emph{pre-parity complex}.
We chose to implement this concept first, as there are many trivial results about preparity complexes that we will later use.
A preparity complex is implemented as the following data:
\begin{code}
C : Type
dim : C $\rightarrow$ nat
plus : C $\rightarrow$ Ensemble C
minus : C $\rightarrow$ Ensemble C
\end{code}
This data is technically different from our description above, but the essential structure is identical.
There is a collection of objects $C$, each member of which has a dimension and two face-sets\footnote{In the actual code the type \lstinline{C} is called \lstinline{carrier}.}.
A few axioms are introduced to ensure that face-sets are finite, non-empty, and disjoint, and that they interact with dimension correctly.
\begin{code}
forall (x y : C), x $\in$ (plus y)  $\rightarrow$ dim y = dim x + 1
forall (x y : C), x $\in$ (minus y) $\rightarrow$ dim y = dim x + 1 
forall (x : C),   Finite (plus x)
forall (x : C),   Finite (minus x)
forall (x : C),   dim x > 0 $\rightarrow$ Inhabited (plus x)
forall (x : C),   dim x > 0 $\rightarrow$ Inhabited (minus x)
forall (x : C),   dim x = 0 $\rightarrow$  plus x == Empty_set
forall (x : C),   dim x = 0 $\rightarrow$ minus x == Empty_set
forall (x : C),   Disjoint (plus x) (minus x)
\end{code}

These are given meaningful names such as \lstinline{plus_Finite}, \lstinline{plus_dim}, and \lstinline{plus_Inhabited}.
Fundamental definitions for sets such as $S_n$ and $S^n$ are also given and some trivial statements are also proved here.
For example,
\begin{code}
Definition sub (R : Ensemble C) (n : nat) : Ensemble C 
 := fun (x : C) => (x $\in$ R $\wedge$ (dim x) = n).

Lemma sub_Union : 
 forall T R n, 
  sub (T $\cup$ R) n == (sub T n) $\cup$ (sub R n).
\end{code}
More complicated definitions like well-formedness are also given and more powerful (though almost trivial) results are also proved here.
For example,
\begin{code}
Definition well_formed (X : Ensemble C) : Prop :=
 (forall (x y : C), x $\in$ X $\wedge$ y $\in$ X
   $\rightarrow$ dim x = O $\rightarrow$ dim y = 0 
   $\rightarrow$ x = y)
 $\wedge$
 (forall (x y : C), x $\in$ X $\wedge$ y $\in$ X
    $\rightarrow$ (forall (n : nat), dim x = S n $\rightarrow$ dim y = S n 
    $\rightarrow$ $\neg$ (perp x y) $\rightarrow$ x = y)).

Lemma well_formed_by_dimension :
 forall X,
  well_formed X <-> forall n, well_formed (sub X n).
\end{code}
All other basic definitions and trivial results are encoded in a similar fashion.
\end{implementation}

We now look at some basic properties of parity complexes.

Given $S \subseteq C$, let $S^\mp$ denote the set of negative faces of elements of $S$ which are not positive faces of any element of $S$, and similarly for $S^\pm$; thus
$$S^\mp = S^- \setminus S^+ \quad \text{and} \quad S^\pm = S^+ \setminus S^- \rlap{\ .}$$
This extends to individual elements by $x^\pm := \{x\}^\pm$ and $x^\mp := \{x\}^\mp$.
These sets capture the notion of \emph{purely positive} and \emph{purely negative} faces of an element $x$ or set $S$.

The following propositions follow from Axioms 1, 2 and 3.

\begin{proposition}[Proposition 1.1]\label{prop_1_1}
For all $x\in C$,
\begin{equation}
{x}^{++} \cap {x}^{--} = {x}^{-+} \cap {x}^{+-} = \emptyset
\end{equation}
\begin{equation}
{x}^{-\mp} = {x}^{+\mp} = {x}^{--} \cap {x}^{+-}
\end{equation}
\begin{equation}
{x}^{-\pm} = {x}^{+\pm} = {x}^{-+} \cap {x}^{++}\rlap{~.}
\end{equation}
\end{proposition}

Proposition~\ref{prop_1_1} contains identities that one would expect from a polytope-like structure and are much like Axiom~\ref{axiom1}.
The meaning is reasonably clear when the various face-sets are highlighted in an example like~\eqref{axiom1_diagram} above.

\begin{proposition}[Proposition 1.2]\label{prop_1_2}
For all $u,v,x\in C$,
$u \lhd v$ and
$v \in {x}^{+}$ imply
\begin{equation}
{u}^{-} \cap {x}^{-+} = \emptyset\rlap{~.}
\end{equation}
\end{proposition}

Proposition~\ref{prop_1_2} indicates that if $u$ branches out from the source of $x$ then a path from $u$ to $v$ can not end in the target of $x$. 
This is a consequence of Axiom~\ref{axiom3}B and has three duals obtained by reversing the roles of $u$ and $v$ and reversing the roles of $x^-$ and $x^+$. 
Proposition~\ref{prop_1_2} and its duals are together equivalent to Axiom~\ref{axiom3}B.

The following observation describes a convenient technical property of well-formed sets.

\begin{observation}[page 322 in~\cite{Street1991}]\label{obs:A}
For all $T,Z \subseteq C$, if
$T \cup Z$ is well-formed and
$T \cap Z = \emptyset$, then
$T \perp Z$.
\end{observation}

We say a set $R \subseteq C$ is \emph{tight} when, for all $u,v\in C$, $u \lhd v$ and $v \in R$ implies ${u}^{-} \cap R^\pm$ is empty.
This condition prevents a path from starting in $R^\pm$ and ending in $R$.
The following two results are required for somewhat technical reasons.

\begin{definition} Suppose that $R,T \subseteq C$.
We say that \emph{$R$ is a segment of $T$}
when for all $x,y,z \in T$, $x,z\in R$ and $x \lhd y \lhd z$ implies $y\in R$.
\end{definition}

\begin{proposition}[Proposition 1.4]\label{prop_1_4}	
For all $R,S \subseteq C$, if
$R$ is tight, $S$ is well-formed, and $R \subseteq S$, then
$R$ is a segment of $S$.
\end{proposition}

\begin{observation}[page 359 in~\cite{Street1994}]\label{obs:B}
For all $x\in C$, ${x}^{+}$ and ${x}^{-}$ are tight.
\end{observation}

This concludes our exposition of \S1.
We have seen that the content of this section could be implemented with very little deviation from the original text.

\begin{remark}
The notion of tightness was introduced in the Corrigenda~\cite{Street1994}.
It appears to be entirely necessary, but we do not understand the full significance of the concept.
\end{remark}

\begin{implementation}
Each axiom and proposition is readily encoded, for example
\begin{code}
Axiom axiom1 :
 forall (x : C),
  (Plus ( plus x)) $\cup$ (Minus (minus x)) == 
  (Plus (minus x)) $\cup$ (Minus ( plus x)).

Lemma Prop_1_2 :
 forall u v x,
  triangle u v $\rightarrow$
   v $\in$ (plus x) $\rightarrow$
    (minus u) $\cap$ (Plus (minus x)) == Empty_set.
\end{code}
We were able to prove each result from basic definitions and axioms. 
This is exactly as described in the original work.
The proof of Proposition~\ref{prop_1_4} makes use of Proposition~\ref{prop_1_1,prop_1_2}.

When we look ahead we find that axioms 1 and 2 are used frequently throughout the material.
Axiom~\ref{axiom3}A is only used to prove that $\lhd_S$ is decidable and that finite non-empty sets $S$ have minimal and maximal elements under $\lhd_S$.
Axiom~\ref{axiom3}B is used only to prove Proposition~\ref{prop_1_1,prop_1_2} and some disjointness conditions in the proof of Lemma~\ref{lemma_3_2}. 

\end{implementation}

\begin{remark}[Adjusting the axioms]
In private conversation Christopher Nguyen pointed out that Axiom~\ref{axiom3}B is only used to prove Proposition~\ref{prop_1_1}, and Proposition~\ref{prop_1_2} and its duals.
We have commented already that Proposition~\ref{prop_1_2} and its duals are equivalent to Axiom~\ref{axiom3}B.
A quick examination of our code then reveals that Proposition~\ref{prop_1_2} is only used to prove that $x^+$ is tight and the disjointness condition described on p327.
We haven't investigated this in any detail, but it might be possible to replace Axiom~\ref{axiom3}B with something slightly weaker (or stronger) but which has the same implications in the relevant proofs.
This is of particular use in light of the fact that Axiom~\ref{axiom3}B is not always preserved under products and joins (see the remark on page 334 of~\cite{Street1991}).
\end{remark}

\section{Movement}  \label{sec:movement}

In Section 2 of~\cite{Street1991} the concept of \emph{movement} is introduced.
It is a concept that is fundamental to describing cells in $n$-categories generated from parity complexes.
For three sets $S,M,P\subseteq C$, we say that \emph{$S$ moves $M$ to $P$}, or $M\xrightarrow{S}P$, when
\begin{equation}
M = (P \cup {S}^{-}) \setminus {S}^{+}
\quad \text{and} \quad
P = (M \cup {S}^{+}) \setminus {S}^{-}\rlap{~.}
\end{equation}
Here are some examples of movement at dimensions 2 and 1:
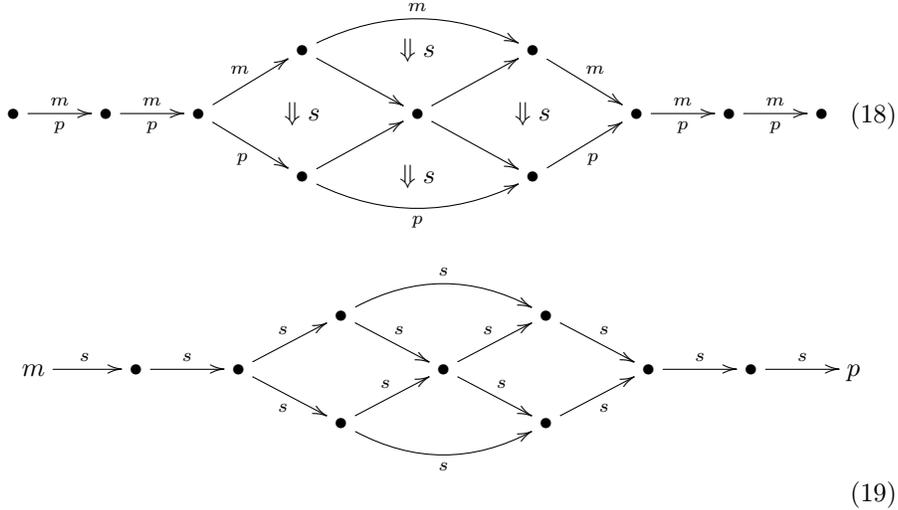
\begin{equation}
\cd[@R-1.5em]{
	& & & \node \ar@/^{1.2em}/[rr]^{m} \ar[dr] & \Downarrow s & \node \ar[dr]^{m} & & &  \\
\node \ar[r]^{m}_{p} & \node \ar[r]^{m}_{p} & \node \ar[ur]^{m} \ar[dr]_{p} & \Downarrow s & \node \ar[ur] \ar[dr] & \Downarrow s & \node \ar[r]^{m}_{p} & \node \ar[r]^{m}_{p} & \node  \\
& &	& \node \ar@/_{1.2em}/[rr]_{p} \ar[ur] & \Downarrow s & \node \ar[ur]_{p} & & &  	 
}
\end{equation}
\begin{equation}
\cd[@!C@R-1.5em]{
	& & & \node \ar@/^{1.2em}/[rr]^{s} \ar[dr]^{s} & & \node \ar[dr]^{s} & & &  \\
m \ar[r]^{s} & \node \ar[r]^{s} & \node \ar[ur]^{s} \ar[dr]_{s} & & \node \ar[ur]^{s} \ar[dr]^{s} & & \node \ar[r]^{s} & \node \ar[r]^{s} & p \\
& &	& \node \ar@/_{1.2em}/[rr]_{s} \ar[ur]^{s} & & \node \ar[ur]_{s} & & &  	 
}
\end{equation}
where lowercase labels $m,p,s$ indicate which set each component belongs to (unlabelled elements do not belong to $M$, $P$, or $S$).
This condition guarantees that the face-sets of $S$, $M$ and $P$ are related in the basic way we would expect of pasting diagrams in $n$-categories.
The movement condition is intended to describe the basic combinatorial shape of cells in our yet-to-be-defined $\omega$-category.
When those cells are defined we will need to add basic finiteness and well-formedness conditions to ensure that various pathogical examples of movement are excluded.

It is helpful to recognise that movement is a condition that applies dimension-by-dimension, that is,
$M\xrightarrow{S}P$ if and only if $M_n\xrightarrow{S_{n+1}}P_n$ for all $n$.
This not only aids in various proofs, but it indicates there is nothing complicated happening across dimensions.

\begin{proposition}[Proposition 2.1]\label{prop_2_1}
	For all $S,M \subseteq C$, there exists $P\subseteq C$ with
	$M\xrightarrow{S}P$
	if and only if
	\begin{equation}
{S}^{\mp} \subseteq M \quad \text{and} \quad \text{$M \cap {S}^{+} = \emptyset$}\rlap{~.}
	\end{equation}
\end{proposition}

Proposition~\ref{prop_2_1} illuminates a fundamental meaning of movement: that $M$ contains the purely negative faces of $S$ and none of the positive faces.
This is illustrated below where elements of $S_{2}^{\mp}$ are indicated by squiggly arrows and those of $S_{2}^{+}$ are indicated by dashed arrows.
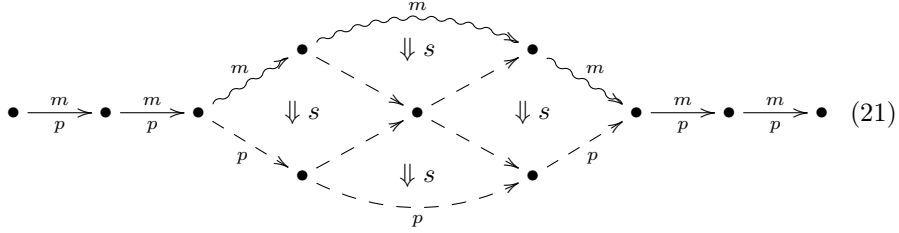
\begin{equation}
\cd[@R-1.5em]{
& & & \node \ar@{~>}@/^{1.2em}/[rr]^{m} \ar@{-->}[dr] & \Downarrow s & \node \ar@{~>}[dr]^{m} & & &  \\
\node \ar[r]^{m}_{p} & \node \ar[r]^{m}_{p} & \node \ar@{~>}[ur]^{m} \ar@{-->}[dr]_{p} & \Downarrow s & \node \ar@{-->}[ur] \ar@{-->}[dr] & \Downarrow s & \node \ar[r]^{m}_{p} & \node \ar[r]^{m}_{p} & \node \\
& &	& \node \ar@{-->}@/_{1.2em}/[rr]_{p} \ar@{-->}[ur] & \Downarrow s & \node \ar@{-->}[ur]_{p} & & &  	 
}
\end{equation}
Observe that $S_{2}^{\mp} \subseteq M_1$ and $M_1 \cap S_{2}^{+} = \emptyset$ as indicated by the proposition.
Proposition~\ref{prop_2_1} has a dual where $M$ and $P$ play opposite roles.

\begin{proposition}[Proposition 2.2]\label{prop_2_2}
Suppose $S,M,P,X,Y \subseteq C$, 
$M\xrightarrow{S}P$ and
$X \subseteq M$ has ${S}^{\mp} \cap X = \emptyset$.
If $Y \cap {S}^{+} = \emptyset$, and
$Y \cap {S}^{-} = \emptyset$, then
$(M{\cup}Y){\cap \neg }X\xrightarrow{S}{(P{\cup}Y){\cap \neg }X}$.
\end{proposition}
Proposition~\ref{prop_2_2} indicates that some elements of $M$ and $P$ can be added or removed without disturbing the movement condition.
The conditions on $X$ and $Y$ indicate that they are disjoint from the faces of $S$ in a suitable way.
Sets $X$ and $Y$ should be thought of as sets that are added to or removed from the movement as below.
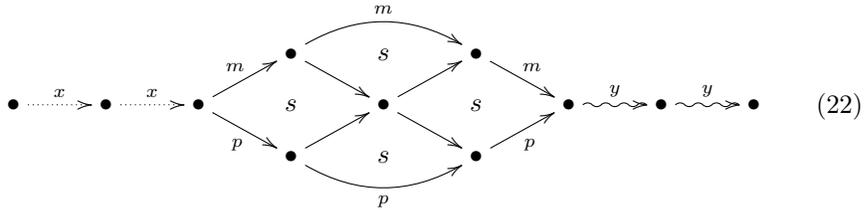
\begin{equation}
\cd[@R-1.5em]{
	& & & \node \ar@/^{1.2em}/[rr]^{m} \ar[dr] & s & \node \ar[dr]^{m} & & &  \\
\node \ar@{.>}[r]^{x} & \node \ar@{.>}[r]^{x} & \node \ar[ur]^{m} \ar[dr]_{p} & s & \node \ar[ur] \ar[dr] & s & \node \ar@{~>}[r]^{y} & \node \ar@{~>}[r]^{y} & \node  \\
& &	& \node \ar@/_{1.2em}/[rr]_{p} \ar[ur] & s & \node \ar[ur]_{p} & & &  	 
}
\end{equation}

\begin{proposition}[Proposition 2.3]\label{prop_2_3}
Suppose $M,P,Q,S,T \subseteq C$ where
$M\xrightarrow{S}P$ and
$P\xrightarrow{T}Q$.
If ${S}^{-} \cap {T}^{+} = \emptyset$ then
$M~\xrightarrow{S \cup T}~Q$.
\end{proposition}

Proposition~\ref{prop_2_3} describes the condition under which movements can be `composed' or `pasted' together.
The following diagram depicts an example. 
Elements of sets $M,S,P,T,Q$ are labelled with the corresponding lower-case letters.
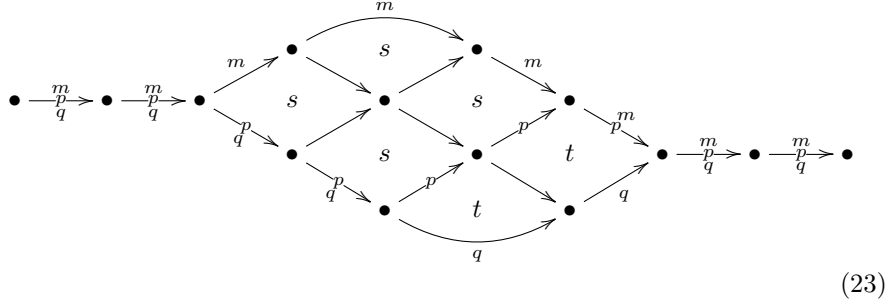
\begin{equation}
\cd[@R-1.5em]{
	& & & \node \ar@/^{1.2em}/[rr]^{m} \ar[dr] & s & \node \ar[dr]^{m} & & &  \\
\node \ar[r]^{m}|{p}_{q} & \node \ar[r]^{m}|{p}_{q} & \node \ar[ur]^{m} \ar[dr]|{p}_{q} & s & \node \ar[ur] \ar[dr] & s & \node \ar[dr]^{m}|{p} & &  \\
& & & \node \ar[ur] \ar[dr]|{p}_{q} & s & \node \ar[ur]|{p} \ar[dr] & t & \node \ar[r]^{m}|{p}_{q} & \node \ar[r]^{m}|{p}_{q} & \node  \\
&  & & 	& \node \ar@/_{1.2em}/[rr]_{q} \ar[ur]|{p} & t & \node \ar[ur]_{q} & & &  	 
}
\end{equation}

\begin{proposition}[Proposition 2.4]\label{prop_2_4}
Suppose $M~\xrightarrow{T \cup Z}~P$ with
${Z}^{\pm} \subseteq P$.
If
$T \perp Z$ then there exists $N$ such that
$M\xrightarrow{T}N\xrightarrow{Z}P$.
\end{proposition}

Proposition~\ref{prop_2_4} describes a condition under which movement can be decomposed.
In particular, if $T \cup Z$ is well-formed then $T \perp Z$ as required in the proposition.

\begin{implementation}
The definition of movement and the propositions above are readily encoded.
For example:
\begin{code}
Definition moves_def (S M P : Ensemble C) : Prop :=
  P == ((M $\cup$ ( Plus S)) $\cap$ (Complement (Minus S)))
  $\wedge$
  M == ((P $\cup$ (Minus S)) $\cap$ (Complement ( Plus S))).

Notation "S 'moves' M 'to' P" := (moves_def S M P) (at level 89).

Lemma Prop_2_3 : forall (S M P T Q : Ensemble C),
  S moves M to P $\rightarrow$ 
  T moves P to Q $\rightarrow$ 
  (Disjoint (Minus S) (Plus T)) $\rightarrow$
    (S $\cup$ T) moves M to Q.
\end{code}
The \lstinline{Notation} command in Coq allows us to use the statement \lstinline{S moves M to P} in place of the somewhat awkward \lstinline{moves_def S M P}.

It did not take long to verify that the proofs in this section proceed precisely as indicated in the original text.

Proposition~\ref{prop_2_1}  is proved by appealing to definitions and basic manipulation of sets.
Proposition~\ref{prop_2_2,prop_2_3,prop_2_4} are proved using Proposition~\ref{prop_2_1} and basic manipulation of sets.
Proposition~\ref{prop_2_1,prop_2_4} have duals that are not displayed here but are required later; they are implemented separately in our code.
It is worth noting that none of these results require axioms 1, 2 or 3.
In our implementation, we prove these results before the axioms are even introduced.

\end{implementation}

This concludes our exposition of \S2.
Again, the content of this section was implemented with very little deviation from the original text. 

\section{The $\omega$-category of a parity complex}  \label{sec:thewcategory}

Having described the basic properties of parity complexes and the more advanced notion of movement, Section 3 of~\cite{Street1991} describes the cells of an $\omega$-category $\O(C)$ associated with any parity complex $C$.

\begin{definition}
A \emph{cell} of a parity complex $C$ is a pair $(M,P)$ of non-empty, well-formed, finite, subsets of $C$ with the property that $M$ and $P$ both move $M$ to $P$.
\end{definition}
If this is interpreted dimension by dimension, we get the following picture at dimension 2, 
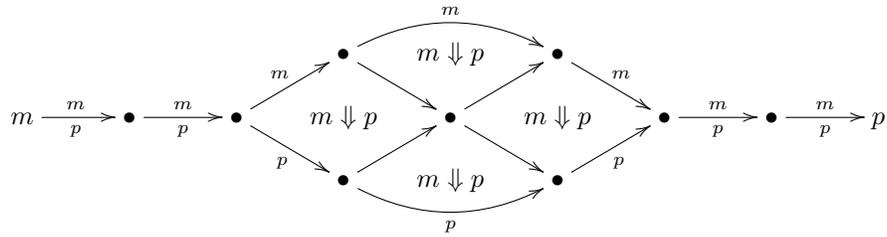
\begin{equation}\label{cell_diagram}
\cd[@!C@C-1.5em@R-1.5em]{
	& & & \node \ar@/^{1.2em}/[rr]^{m} \ar[dr] & m \Downarrow p & \node \ar[dr]^{m} & & &  \\
m \ar[r]^{m}_{p} & \node \ar[r]^{m}_{p} & \node \ar[ur]^{m} \ar[dr]_{p} & m \Downarrow p & \node \ar[ur] \ar[dr] & m \Downarrow p & \node \ar[r]^{m}_{p} & \node \ar[r]^{m}_{p} & p  \\
& &	& \node \ar@/_{1.2em}/[rr]_{p} \ar[ur] & m \Downarrow p & \node \ar[ur]_{p} & & &  	 
}
\end{equation}
where lowercase labels $m,p,s$ indicate which set the elements belong to. 
Notice that $M_1$ and $P_1$ are neither equal nor disjoint, but each move $M_0$ to $P_0$.
Notice also that $M_2=P_2$.
This kind of behaviour is uniform through all dimensions.
Notice also that, aside from the movement condition, we only require that $M$ and $P$ be non-empty, well-formed and finite. 
Call $(M,P)$ an \emph{$n$-cell} when $M \cup P$ is $n$-dimensional.
In this case we have $M_n = P_n$ as above.

\begin{definition}
The \emph{$n$-source} and \emph{$n$-target} of a pair of sets $(M, P)$ are defined by
\begin{equation}
s_n(M,P) = ({M}^{n-1} \cup {M}_{n} , {P}^{n-1} \cup {M}_{n})
\end{equation}
and
\begin{equation}
t_n(M,P) = ({M}^{n-1} \cup {P}_{n} , {P}^{n-1} \cup {P}_{n})\rlap{\rlap{~.}}
\end{equation}
\end{definition}
If $(M,P)$ is a cell we can show that $s_n(M,P)$ and $t_n(M,P)$ are also cells, and that they are $n$-dimensional. 
Notice that $(M,P)$, $s_n(M,P)$ and $t_n(M,P)$ contain exactly the same elements in dimension $n-1$ and below.
We encourage the reader to consider the 1-source and 1-target of the cell depicted in~\eqref{cell_diagram}.

\begin{definition}
A pair of cells $(M, P)$, $(N, Q)$ are \emph{$n$-composable} when
\begin{equation}\label{source_equal_target}
t_n(M,P) = s_n(N,Q)\rlap{\ ,} 
\end{equation}
in which case their \emph{$n$-composite} is
\begin{equation}
(N, Q) *_n (M, P) := (M \cup (N \cap \neg  {N}_{n}),\ (P \cap \neg {P}_{n}) \cup Q)\rlap{~.}
\end{equation}
\end{definition}
Notice that~\eqref{source_equal_target} implies that the two cells agree from dimensions 0 to $n-1$ and that $P_n = N_n$ at dimension $n$.
The resulting composite is almost exactly the pair-wise union of $(M,P)$ and $(N,Q)$; the set-difference ensures correct behaviour at dimension $n$.
It is not surprising that some form of set-difference is required since most forms of composition will forget the point of contact: $ A \to B \to C$ composes to $A \to C$.

For any parity complex $C$, let $\O(C)$ be the set of cells of $C$.
We will see later (Theorem~\ref{theorem_3_6}) that $\O(C)$ is an $\omega$-category.
Before this can be achieved, we need to establish some basic properties of cells.

\begin{definition}
A set $S \subseteq C$ is \emph{receptive} when for all $x\in C$,
$$ 	\text{if} \quad
	{x}^{-+} \cap {x}^{++} \subseteq S
	\quad \text{and} \quad
	S \cap {x}^{--} = \emptyset
	\quad\text{then}\quad
	S \cap {x}^{+-} = \emptyset$$
and
$$ 	\text{if} \quad
	{x}^{+-} \cap {x}^{--} \subseteq S
	\quad\text{and}\quad
	S \cap {x}^{++} = \emptyset
	\quad\text{then}\quad
	S \cap {x}^{-+} = \emptyset\rlap{\rlap{~.}}$$
\end{definition}
A cell is receptive when it is receptive at every dimension.

\begin{remark}
The notion of receptivity is somehow important, we find later that all cells are receptive and it is a necessary condition for some central results.
It appears to be entirely necessary, but we do not have an intuitive understanding of its meaning.
\end{remark}

\begin{lemma}[Lemma 3.1]\label{lemma_3_1}
For all $M,P\subseteq C$, $x \in C$, if
$M\xrightarrow{{x}^{+}}P$ and
$M$ is receptive then
$M\xrightarrow{{x}^{-}}P$.
\end{lemma}
Lemma~\ref{lemma_3_1} is proved using definitions, basic manipulation of sets and Proposition~\ref{prop_2_1,prop_1_1}.
It has a dual which we implement in our code.
We will find later that since all cells are receptive, it is not hard to find receptive subsets $M$ of $C$.
In fact, it is a bit difficult to illustrate why receptivity is even required because the most obvious examples of $x,M,P$ satisfying the movement condition above are also part of a cell structure.

\begin{lemma}[Lemma 3.2]\label{lemma_3_2}
Suppose $m,n \in \mathbb{N}$, all cells are receptive and $(M, P)$ is an $n$-cell.
Suppose also that $X \subseteq C_{n+1}$, $\abs{X} = m$ and $X$ is well-formed with
${X}^{\pm} \subseteq {M}_{n}$.
Put $Y = ({M}_{n} \cup {X}^{-}) \cap \neg  {X}^{+}$, then:
\begin{enumerate}[label=\Alph*., itemindent=3.5em, labelsep=1.5em]
\setcounter{enumi}{1}
\item $({M}^{n-1} \cup Y, {P}^{n-1} \cup Y)$ is a cell and
and ${X}^{-} \cap {M}_{n} = \emptyset$.
\item  $({M}^{n-1} \cup Y \cup X,
P \cup X)$ is a cell.
\end{enumerate}
\end{lemma}
Lemma~\ref{lemma_3_2} originally contained a part A which was removed in~\cite{Street1994}.
Lemma~\ref{lemma_3_2}C indicates that, if $X$ is a well-formed set of dimension $n+1$, $(M,P)$ is an $n$-cell, and $X$ abuts $(M,P)$ in the sense that $X^\pm \subseteq M_n$, then we can form an $(n+1)$-cell whose top-dimension elements are those of $X$ and whose target is $(M,P)$. 
The source of this cell has $Y$ at its top dimension.
The following diagram is labelled to illustrate this scenario.
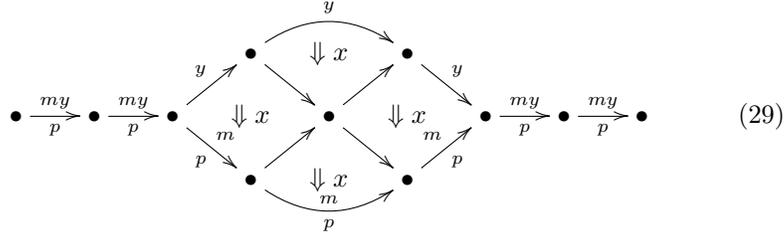
\begin{equation}
\cd[@!C@C-1.5em@R-1.5em]{
	& & & 
	\node \ar@/^{1.2em}/[rr]^{y} \ar[dr] & \Downarrow x & \node \ar[dr]^{y} & & &  \\
\node \ar[r]^{my}_{p} & \node \ar[r]^{my}_{p} & \node \ar[ur]^{y} \ar[dr]^{m}_{p} & \Downarrow x & \node \ar[ur] \ar[dr] & \Downarrow x & \node \ar[r]^{my}_{p} & \node \ar[r]^{my}_{p} & \node  \\
& &	& \node \ar@/_{1.2em}/[rr]^{m}_{p} \ar[ur] & \Downarrow x & \node \ar[ur]^{m}_{p} & & &
}
\end{equation}
There is a dual lemma obtained by reversing the direction of $X$ in the diagram above.

This kind of result does not seem unusual, but it is surprisingly hard to prove (see Implementation~\ref{implementation_3_2} below).
The proof itself is done in three steps. 
To quickly summarise:
\begin{enumerate}
\item Lemma~\ref{lemma_3_2}B implies Lemma~\ref{lemma_3_2}C.
The proof is somewhat direct and proceeds as indicated in the original paper.
\item Lemma~\ref{lemma_3_2}B with $m=1$ implies Lemma~\ref{lemma_3_2}B in general.
This is done by induction on $m$ and follows from basic definitions and axioms.
\item Lemma~\ref{lemma_3_2}B holds for $m=1$.
This is done by induction on $n$ and the argument relies on Proposition~\ref{prop_3_3}.
The construction works as indicated, though it is not a short argument. 
There are particular disjointness conditions that must be established (p327 of~\cite{Street1991}) and require their own special argument.
\end{enumerate}

\begin{proposition}[Proposition 3.3]\label{prop_3_3}
For all $n \in \mathbb{N}$, all $n$-cells in $C$ are receptive.
\end{proposition}
This is a somewhat technical result, it is not immediately clear to us how the notion of receptivity fits naturally into the combinatorics.
The proof of this result relies on Lemma~\ref{lemma_3_2}B.

\begin{theorem}[Theorem 3.6]\label{theorem_3_6}
If $C$ is any parity complex then $\O(C)$ is an $\omega$-category. 
Furthermore, if $(M,P)$, $(N,Q)$ are $n$-composable cells\footnote{$t_n(M,P) = s_n(N,Q)$} then $ (M_k \cup P_k)^- \cap (N_k \cup Q_k)^+ = \emptyset$
for all $k > n$.
\end{theorem}

Theorem~\ref{theorem_3_6} is a central result in~\cite{Street1991} since it achieves one of the main goals of the paper.
In order to implement Theorem~\ref{theorem_3_6} we would need to implement a notion of $\omega$-category which is not trivial.
Since there is little question that this result holds, and it is not required to prove Theorem~\ref{theorem_4_1}, we have chosen not to implement it.
We similarly omit Propositions 3.4 and 3.5 which are preliminary results leading up to Theorem~\ref{theorem_3_6}.

%

\begin{remark}
Perceptive readers will have noticed that Lemma~\ref{lemma_3_2}B and Proposition~\ref{prop_3_3} seem to logically rely on one another.
At first glance this appears to be a circular argument and therefore unsound.
However, if we look closely we can see that each result proceeds by induction and that the two proofs can be woven together to produce a proof of both results simultaneously. 
Proposition~\ref{prop_3_3} is restated as: \emph{for all $n$, every $n$-dimensional cell $(M,P)$ is receptive}.
The two results are proved by mutual induction on $n$, the dimension of $(M,P)$.
Included in that argument is an induction on $m = \abs{X}$.
The following statements hold and are enough to show that both results hold for all $n$ and $m$.
\begin{enumerate}[label=\roman*.]
\item Lemma~\ref{lemma_3_2}B holds when $m=1$ and $n=0$.
\item For a fixed $n$, if Lemma~\ref{lemma_3_2}B holds when $m=1$ then it holds when $m>1$ (by induction on $m$).
\item Proposition~\ref{prop_3_3} holds when $n=0$.
\item If Lemma~\ref{lemma_3_2}B and Proposition~\ref{prop_3_3} hold for $n=k$, then Lemma~\ref{lemma_3_2}B holds for $n=k+1$ and $m=1$.
\item If Lemma~\ref{lemma_3_2}B holds for $n=k+1$ and Proposition~\ref{prop_3_3} holds for $n=k$, then Proposition~\ref{prop_3_3} holds for $n=k+1$.
\end{enumerate}
This understanding is not explicit in~\cite{Street1991}.
\end{remark}

\begin{implementation}\label{implementation_3_2}
As in earlier sections, the definitions and statement of results are readily encoded.
The main difficulty arises in encoding the proofs.

The proofs of Lemma~\ref{lemma_3_2} and Proposition~\ref{prop_3_3} are by far the most difficult part of the entire project and consumed most of our programming effort. 
Consider the components of the proof given above. 
Each of the components follow the argument provided by Street in his paper. 
However the disjointness condition in (iv) has a dual, and (i), (ii), and (iv) each have duals. 
Finally, we needed to uncover the logical dependence that allows us to weave these things together to produce a non-cyclic argument. 

It is worth noting that the original proof of Proposition~\ref{prop_3_3} uses an argument about skeletons of parity complexes (treating separate parity complexes as objects of the argument).
We have translated the argument so that it is internal to any given parity complex.
The combinatorial logic of our argument is exactly the same as Street's, we have only adjusted the setting slightly.
\end{implementation}

\begin{remark}[Understanding receptivity and tightness]
This section of~\cite{Street1991} is harder to follow than the others and the proofs here are not straight-forward.
In particular, the notions of tightness and receptivity are both a bit opaque and Lemma~\ref{lemma_3_2} is very hard to prove.
This provides some motivation to closely re-examine the result and see whether alternative arguments might be made to prove it.
\end{remark}

This concludes our exposition of \S3. In formalizing this section we have seen some of the most complicated arguments made in the original text and seen that Lemma~\ref{lemma_3_2} (Lemma 3.2) is not easy to prove.

\section{Freeness of the $\omega$-category} \label{sec:freeness}

Having built the $\omega$-category $\O(C)$ from a parity complex $C$, we now prove that it is generated from atoms. The following content comes from \S4.

In any parity complex $C$ we expect that any individual element $x$ of dimension $p$ is the top element of some cell whose lower-dimensional structure can be computed by examining the face-sets of $x$ and recursively taking face-sets of face-sets.
This is made explicit in the following definition.

\begin{definition}
For each $x \in C_p$, two subsets $\mu(x),\pi(x) \subseteq C^p$ are defined inductively as follows
$$ \mu(x)_p = \{x\} \quad \text{and} \quad \mu(x)_{k-1} = \mu(x)_{k}^\mp, \quad 1 \leq k \leq p$$
$$ \pi(x)_p = \{x\} \quad \text{and} \quad \pi(x)_{k-1} = \pi(x)_{k}^\pm, \quad 1 \leq k \leq p$$
The pair $(\mu(x), \pi(x))$ is denoted by $\braket{x}$.
\end{definition}

Take the following diagram for example.
If $x\in C$ has dimension 2 and has boundary as illustrated in below then $\braket{x} = (\{x,p,q,r,a\},\{x,s,t,e\}) $
\begin{equation}\label{mu_pi_example}
\cd[]{
 & 
 c \ar[r]^{q} & 
 b \ar[dr]^{r} & \\ 
 a \ar[ur]^{p} \ar[r]_{s} \dtwocell{urrr}{x} & d \ar[rr]_{t} & & e
 }
\end{equation}
A priori, we have no guarantee that such a pair is actually a cell.

\begin{definition}
An element $x \in C_p$ is called \emph{relevant} when $\braket{x}$ is a cell. 
This amounts to saying that $\mu(x)_n$ and $\pi(x)_n$ are well-formed for $0 \leq n < p-1$, and
$$ \mu(x)_{n-1} = \pi(x)_n^\mp, \qquad \pi(x)_{n-1} = \mu(x)_n^\pm $$
for $0 < n < p-1$.
Call a cell $(M,P)$ an \emph{atom} when it is equal to $\braket{x}$ for some $x \in C$. 
In that case we say that $(M,P)$ is \emph{atomic}.
\end{definition}

In all of our main examples, every $\braket{x}$ is a cell (all elements are relevant).

\begin{theorem}[Theorem 4.1: excision of extremals]\label{theorem_4_1}
Suppose that $\mu(x)$ is tight for all $x \in C$.
Suppose $(M , P)$ is an $n$-cell and $u \in M_n$ ($= P_n$) is such that $(M, P) \not= \braket{u}$~\footnote{Alternatively, let $(M,P)$ be a non-atomic $n$-cell.}.
Then $(M, P)$ can be decomposed as
\begin{equation}
	(M, P) = (N, Q) *_m (L, R)
	\end{equation}
where $m<n$, and $(N,Q)$ and $(L,R)$ are $n$-cells of dimension greater than $m$.
\end{theorem}
This is another central result of the paper.
If this algorithm is applied recursively then it shows how to present an arbitrary $n$-cell as a composite of atoms. 
Thus $\O(C)$ is not only an $\omega$-category, but it is generated from its atoms

The algorithm takes an $n$-cell $(M, P)$ and runs as follows.
\begin{enumerate}
\item Find the largest $m<n$ with $(M_{m+1}, P_{m+1}) \not= (\mu(u)_{m+1}, \pi(u)_{m+1})$.
This amounts to discovering the highest dimension at which the criterion for being atomic does not hold\footnote{Alternatively, find the largest $m<n$ with $M_{m+1} \cap P_{m+1} \not= \emptyset$.}.
In this case, there exists $w \in M_{m+1} \cap P_{m+1}$.
\item We want to decompose our cell by pulling off a cell of dimension $m+1$.
Let $x$ be a minimal element of $M_{m+1}$ less than $w$, and let $y$ be a maximal element of $M_{m+1}$ greater than $w$.
\item At least one of $x$ or $y$ must belong to $M_{m+1} \cap P_{m+1}$. 
This relies on the fact that $\mu(u)_{m+1}$ is a segment of $M_{m+1}$, which itself relies on $\mu(u)_{m+1}$ being tight.
\item If $x \in M_{m+1} \cap P_{m+1}$ then we get a decomposition of $(M,P)$ as
\begin{equation}
N = M^{m} \cup \{x\}
\qquad
Q = P^{m-1} \cup ((M_m\cup {x}^{+}) \cap \neg  {x}^{-}) \cup \{x\}
\end{equation}
\begin{equation} 
L = ((M \cap \neg \{x\}) \cup {x}^{+}) \cap \neg {x}^{-}
\qquad
R = P \cap \neg  \{x\}
\end{equation}
Notice that $(N, Q)$ is an $(m+1)$-cell whose single element at top dimension is $x$,  and $(L, R)$ is the $n$-cell obtained by cutting $x$ out of $(M,P)$.
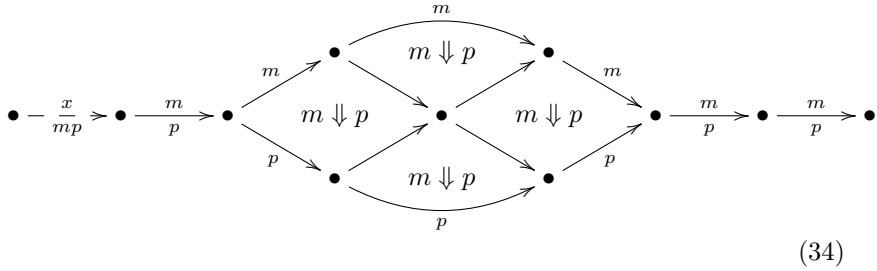
\begin{equation}
\cd[@!C@C-1.5em@R-1.5em]{
	& & & \node \ar@/^{1.2em}/[rr]^{m} \ar[dr] & m \Downarrow p & \node \ar[dr]^{m} & & &  \\ 
\node \ar@{-->}[r]^{x}_{mp} & \node \ar[r]^{m}_{p} & \node \ar[ur]^{m} \ar[dr]_{p} & m \Downarrow p & \node \ar[ur] \ar[dr] & m \Downarrow p & \node \ar[r]^{m}_{p} & \node \ar[r]^{m}_{p} & \node  \\
& &	& \node \ar@/_{1.2em}/[rr]_{p} \ar[ur] & m \Downarrow p & \node \ar[ur]_{p} & & &  	 
}
\end{equation}

\item If $y \in M_{m+1} \cap P_{m+1}$ then we get a decomposition of $(M,P)$ as
\begin{equation}
N = M \cap \neg  \{y\}
\qquad
Q = ((P \cap \neg \{y\}) \cup {y}^{-}) \cap \neg {y}^{+}
\end{equation}
\begin{equation} 
L = M^{m-1} \cup ((P_m \cup {y}^{-}) \cap \neg {y}^{+}) \cup \{y\}
\qquad
R = P^{m} \cup \{y\}
\end{equation}
This is dual to the case for $x$.
Notice that $(L, R)$ is an $(m+1)$-cell whose single element at top dimension is $y$, and $(N,Q)$ is the cell obtained by cutting $y$ out of $(M,P)$. 
\end{enumerate}

The two hardest parts of this algorithm are parts (3) and (4). 
In part (3) we must show that either $x$ or $y$ belong to $M_{m+1} \cap P_{m+1}$.
This relies on the fact that $\mu(x)_{m+1}$ is a segment of $M_{m+1}$, but this follows from Proposition~\ref{prop_2_4} and the assumption that each $\mu(x)$ is tight.
In part (4) we need to show that $(N,Q)$ and $(L,R)$ are well-defined cells. 
The various conditions of finiteness and well-formedness follow quite directly. 
The difficulty comes in showing that the movement conditions hold. 
We investigate the cells dimension by dimension and find that the movement conditions can be proved using~Proposition~\ref{prop_2_4} and~Lemma~\ref{lemma_3_2}.

How do we know that this algorithm terminates? 
The original text defines the rank of an $n$-cell $(M,P)$ to be the cardinality of $M \cup P$.
The algorithm produces two cells of smaller rank, so therefore must terminate.
It is also possible to define the rank by
\begin{equation}
\mathrm{rank}(M,P) = \sum_{k=0}^{n} \vert M_k \cap P_k \rvert
\end{equation}
In this case every $n$-cell has a rank of at least 1 since $M_n \cap P_n$ is non-empty.
A cell of rank $1$ must be atomic.
A cell of rank $k>1$ can be decomposed using excision of extremals into two cells whose individual ranks are less than or equal to $k-1$.
Again, this is sufficient to guarantee termination.

\begin{implementation}
As already indicated, Theorem~\ref{theorem_4_1} is readily proved using the argument given above.
\end{implementation}

\begin{remark}\label{rem:free}
In order to show that $\O(C)$ is \emph{freely} generated from its atoms we must show that there are no equalities among composites of cells that are not a consequence of the $\omega$-category axioms. 
This is achieved in Street's Theorem 4.2 but has not been reproduced here and we have not included it in our formalisation.
\end{remark}

\begin{remark}\label{loose_end_comment1}
Many of these theorems and lemmas come with a condition concerning tightness and receptivity of various sets. 
We can show that these conditions are satisfied by appealing to various other results.
At the end of the day  there may be some confusion about which conditions are ultimately required.
To summarise, if a parity complex $C$ has the property that $\mu(x)$ is tight for every $x \in C$, then all of the theorems up to this point will hold.

At this stage, we have not shown that every $\braket{x}$ is a cell.
In fact, we have no guarantee that any cells exist at all.
This is something of a loose end, it is accounted for in the following section.
\end{remark}

This concludes our exposition of \S4. 
As already mentioned, the central theorem Theorem~\ref{theorem_4_1} (Theorem 4.1) is readily proved using the content of previous sections.

\section{Product and join} \label{sec:product}

Though we have not formalised Section 5 of~\cite{Street1991}, some of the content is relevant to our investigation. 
In \S5 Street describes, for any two parity complexes $C$ and $D$, their \emph{product} $C \times D$ and their \emph{join} $C \bullet D$.
This section also describes two kinds of duals for parity complexes obtained by reversing the roles of $(-)^+$ and $(-)^-$ in all dimensions or in odd dimensions only.
This is of particular interest since the diagrams involved in descent are products of globes with simplexes; this is explored in \S6.

This section also addresses some issues that are, as yet, unresolved.
First, we don't know that any elements are relevant (consequently we don't know if any cells exist at all). 
Second, Lemma~\ref{lemma_3_2} relies on the fact that all $\mu(x)$ are tight, and this was never established.

Consider the following stronger forms of Axioms 1 and 2.
\begin{align*}
\text{For all $x$,} & \\
(R1) \qquad & \mu(x)^- \cup \pi(x)^+ = \mu(x)^+ \cup \pi(x)^- \text{ and } \\
& \mu(x)^- \cap \pi(x)^+ = \mu(x)^+ \cap \pi(x)^- = \emptyset \\
(R2) \qquad & \mu(x) \text{ and } \pi(x)  \text{ are well formed.}
\end{align*}
These axioms hold for $\Delta$, $\glob$, and $\cube$.

\begin{remark}
If a parity complex $C$ satisfies these axioms then every $\braket{x}$ is a cell (every $x$ is relevant). 
Thus all elements of $\Delta$, $\glob$, and $\cube$ are relevant.
\end{remark}

This solves the first problem in the primary examples of interest. 
Now what about tightness of $\mu(x)$?

In a parity complex $C$, write $x \prec y$ when either $y \in x^+$ or $x \in y^-$. 
Let $\blacktriangleleft$ denote the reflexive transitive closure of the relation $\prec$.
Notice that $x < y$ means there exists $z \in x^+ \cap y^-$, so this implies $x \prec y$.
Hence, $x \lhd y$ implies $x \blacktriangleleft y$. 
Where the first ordering $\lhd$ describes paths in a fixed dimension, this new ordering $\blacktriangleleft$ zig-zags up and down with $\prec$ to describe paths across multiple dimensions.
We introduce the following as an optional axiom.
\begin{equation}
(AS) \qquad \text{$\blacktriangleleft$ is anti-symmetric.}
\end{equation}
This axiom holds in $\Delta$, $\glob$, and $\cube$ where $\blacktriangleleft$ is also total.

\begin{proposition}[Proposition 5.2]\label{prop_5_2}
If each $x$ is relevant and (AS) holds then each $\mu(x)$ is tight.
Thus, every $\mu(x)$ in $\Delta$, $\glob$, and $\cube$ are tight.
\end{proposition}

This solves the second problem in the primary examples. Unfortunately we find that there are examples of parity complexes where $\blacktriangleleft$ is not antisymmetric (p337 of~\cite{Street1991}).
These are small pasting diagrams that are explicitly illustrated in the article and are quite elementary.
These stronger conditions do not appear to be unreasonable, but it is not clear what classes of examples are exluded as a consequence.

\begin{remark}\label{loose_end_comment2}
To summarise again, if a parity complex $C$ satisfies (R1), (R2) and (AS) then every theorem and proposition covered here will work.
In particular, every theorem and proposition holds for the parity complexes $\Delta$, $\glob$, and $\cube$.
\end{remark}

\begin{remark}
Theorem~\ref{theorem_4_1} relies on the fact that each $\mu(x)$ is tight and therefore a segment in the required place. 
This is readily proved when (R1), (R2) and (AS) are used, but we would like to use the excision of extremals more generally if possible.
So we ask, is the tightness condition strictly necessary? 
Or, is there another way to ensure that $\mu(x)$ is a segment in that proof?
We do not yet know.
\end{remark}

\begin{remark}
There seems to be a fundamental relationship between parity complexes and `directed graphs of multiple dimension'.
Note that this notion of higher-dimensional graph would not be the same as an $n$-graph since each component of an $n$-graph has a single source and single target rather than a source set and a target set.
Some of the axioms for parity complexes are just those of this `graph' structure and others restrict us to graphs of a certain kind.
Axioms (R1), (R2) and (AS) place further restrictions.
Since we are mainly interested in examples that satisfy all of these conditions, we do not need to worry too much about this narrowing of our focus.
More generally though, it would be good to know which of these conditions are associated with the graph structure of parity complexes, and which of the conditions allow for the (free) $\omega$-category construction.
This could be the focus of some future research.
\end{remark}

This concludes our exposition of \S5. 
We will not discuss Section 6 of~\cite{Street1991} and our exposition of the original text ends here.

\section{Some lessons in coded mathematics} \label{sec:lessons}

\subsection{Duals} 
We were often forced to prove dual results where $x^+$ and $x^-$ were interchanged, or where the direction of a movement $M\xrightarrow{S}P$ was reversed. 
In these cases we were forced to explicitly restate and reprove the result, even though the underlying logic had not changed whatsoever.
It would have been better if, from the beginning, we had encoded plus and minus as duals to each other, then the theorems would dualise automatically. 
One way to do this is to define \lstinline{faceset : bool $\rightarrow$ C $\rightarrow$ Ensemble C } and then set \lstinline{minus := faceset false} and \lstinline{plus := faceset true}. 
From this starting point it should be easy to combine dual results into one.

\subsection{Notation} 
Coq has a \lstinline{Notation} facility which allows the user to introduce custom notation for specific expressions. 
We used this to make set operations easier to read and write.
For example, an expression such as \lstinline{Union A B} is displayed as \lstinline{A $\cup$ B}, and similarly for intersection, inclusion, etc.
This made our code much easier to read.

\subsection{Tactics} 
Coq has a tactic language which allows for partial automation of proofs. 
The language allows the user to describe simple proof strategies that can be automatically applied when little innovative thinking is required.
A particular built-in tactic called \lstinline{intuition} will automatically deal with simple proofs that require only knowledge of first-order logic.
We used the tactic language to describe a proof tactic called \lstinline{basic} that automatically applied further logical steps such as $(x \in A \cap B) \to (x \in A \wedge x \in B)$. 
In many cases this vastly simplified proofs by applying \lstinline{repeat (basic; intuition)} to automatically prove some trivial facts.

More expert use of this system would surely result in more elegant and readable proofs.

\subsection{Axiom of extensionality for sets} 
We chose to remove the axiom of extensionality because we wanted to deal with sets in a completely constructive fashion.
This was a choice of style. 
In many ways, retaining the axiom would not have weakened our encoding and we would not have needed to implement setoid rewrite for ensembles.

\subsection{Compiling the excision of extremals algorithm}
Our choice to implement sets using ensembles has made it impossible to directly compile an executable version of the excision of extremals.
This is unfortunate: we have proved that such an algorithm can run but we can't actually compile or run it without further coding.
The mathematical significance of our work is not undermined, but better planning would have yielded executable code as a pleasant side-effect.

\section{Conclusion} \label{sec:conclusion}
We have formalised Ross Street's \emph{Parity Complexes} up to the excision of extremals algorithm in \S4. 
In particular, \S1 and \S2 together with~Theorem~\ref{theorem_4_1} are proved as indicated in the original text.
Section 3 is also formalised with the same essential arguments as~\cite{Street1991}, but with many additional dual theorems, and a technical but meaningful change to the logical flow of Lemma~\ref{lemma_3_2} and Proposition~\ref{prop_3_3}.

We have indicated where the material is most effective at capturing the difficult combinatorics, and where future work might make improvements. 
We have explicitly outlined the logical dependence of the central results.

\printbibliography

\end{document}